\documentstyle[12pt,amstex,amssymb]{article}

\def\be{\begin{equation}}
\def\ee{\end{equation}}
\def\h{\hbar}

\def\dim{\operatorname{dim}}

\def\ev{\operatorname{ev}}
\def\ft{\operatorname{ft}}

\newcommand{\CC}{{\Bbb C}}

\newcommand{\QQ}{{\Bbb Q}}

\newcommand{\lan}{\langle}
\newcommand{\ran}{\rangle}

\renewcommand{\a}{\alpha}
\renewcommand{\b}{\beta}
\renewcommand{\c}{\gamma}
\renewcommand{\d}{\delta}
\newcommand{\e}{\varepsilon}

\renewcommand{\l}{\lambda}

\newcommand{\z}{\zeta}

\renewcommand{\S}{\Sigma}

\newcommand{\calJ}{{\cal J}}

\newcommand{\calo}{{\cal O}}

\newcommand{\M}{\bar{\cal M}}

\title{The Mirror Formula for Quintic Threefolds}

\author{Alexander Givental \thanks{ Research supported by NSF grant
DMS-9704774} \\ UC Berkeley}

\date{July 14, 1998}

\begin{document}

\maketitle

The formula of Candelas, de la Ossa, Green and Parkes \cite{COGP} 
expressing the virtual numbers $n_d,\ d=1,2,3,...$, of degree $d$ 
holomorphic spheres on quintic threefolds in $\CC P^4$ in terms of
series solutions to the linear differential equation
\be   (q\frac{d}{dq})^4\ \ I = 5q\ 
(5q\frac{d}{dq}+1)(5q\frac{d}{dq}+2)(5q\frac{d}{dq}+3)(5q\frac{d}{dq}+4) 
\ \ I  \label{ODE} \end{equation} 
has been intriguing algebraic and symplectic geometers since the beginning
of the decade. The first proof of this formula was given two years ago 
in the extensive paper \cite{Gi1} among a number of other theorems on 
equivariant Gromov -- Witten theory. Several authors
managed to adjust the approach of that paper to complete intersections in 
homogeneous 
K\"ahler spaces \cite{K1, BCKV, K2}, in toric manifolds \cite{Gi3, Gi4} and to 
symmetric products of Riemann surfaces \cite{BT}.  

We present here a shortcut to our original proof in the case 
of quintic threefolds.
Several variants of the proof can be 
found in \cite{Gi1, Gi3, Gi4, LLY, CP, P} as particular cases of more general 
theorems. Yet it seems useful to illustrate all ingredients of the proof 
in the simplest nontrivial example.

We will assume that the reader is familiar with generalities on orbifolds
and orbibundles, equivariant cohomology and localization formulas,  
Kontsevich's moduli spaces of stable maps  \cite{Kn, BM} and with the 
formulation of the conjecture. We will concentrate therefore
only on the issues relevant for the proof of the mirror formulas. 

In the last section named {\em Updates} we outline several modifications of the
proof available at the moment and provide corresponding references.

\bigskip

{\bf The linear sigma-model.} 
A quintic threefold in $X=\CC P^4$ is given by a generic degree $5$ 
homogeneous equation $Q(x_1,...,x_5)=0$. A degree $d$ parametrized
curve $\CC P^1\to X$ is described by $5$ relatively prime degree $d$ binary
forms which we will represent by the polynomials $x_1(\z),...,x_5(\z)$ of
degree $\leq d$ in the affine coordinate $\z$ on $\CC P^1$. The curve is
situated on the quintic if and only if the degree $\leq 5d$ polynomial
$Q(x_1(\z),...,x_5(\z))$ vanishes identically. This identity yields $5d+1$
equations of degree $5$ in the projective
space $LX_d=\CC P^{5d+4}$ of all (not necessarily relatively prime) $5$-tuples
of degree $\leq d$ polynomials. Attempting to count degree $d$ spheres on
the quintics by means of intersection theory in $LX_d$ we will arrive to the
answer $n_d=5^{5d+1}$ which is meaningful, as it has been explained 
by D. Morrison and R. Plesser in \cite{MP}, but wrong.
   
Our approach to the quintic formula begins with the following observation
\cite{Gi2}. The variety given in $LX_d$ by the equations $Q(x(\z))\equiv 0$
is invariant with respect to the M\"obius transformation group of $\CC P^1$,
and one can employ {\em equivariant} intersection theory for curve counting.
The maximal torus $S^1$ in the M\"obius group acts via $\z\mapsto 
\z \exp{i\phi}$
on the space $LX_d$. The cohomology algebra $H^*_{S^1}(LX_d)$ is generated 
by the equivariant Chern class $-p$ of the Hopf bundle over $LX_d$. 
The polynomial $Q$ defines an invariant section 
$Q(x(\z))$ of an equivariant $5d+1$-dimensional bundle $LV_d$ over $LX_d$. 
The following formal series encodes complete equivariant topological 
information about equivariant Euler classes of these bundles:
\[  L(q,z):=\sum_{d=0}^{\infty} q^d \int_{LX_d} e^{pz} Euler(LV_d)\ . \]
The integral in this formula means evaluation of an equivariant cohomology
class on the invariant fundamental class of the manifold and takes values
in the coefficient algebra $H^*(BS^1)$ of the equivariant cohomology theory.
Thus the coefficients of the $(q,z)$-series $L$ are polynomials in one variable
which we prefer to denote $\h $. With this notation $Euler(LV_d)=
5p (5p-\h)(5p-2\h)...(5p-5d\h)$.

\medskip

{\bf Theorem A} (\cite{Gi2}). {\em (a)
\[  L(q,z)=\lan I(qe^{\h z},\h^{-1}), I(q,-\h^{-1}) \ran \ ,  \]
where 
\[ \lan \phi ,\psi \ran := \int_{\CC P^4}\phi \ \psi \ Euler (\calo (5)) =
\frac{1}{2\pi i} \oint \phi(P)\psi(P) \frac{5P dP}{P^5} \]
is the intersection pairing in the even cohomology algebra 
$\QQ [P]/(P^4)$ of the quintic, and 
\[ I(q,\h^{-1}):=e^{P\ln q/\h } \sum_{d=0}^{\infty } q^d
\frac{(5P+\h)(5P+2\h)...(5P+5d\h)}{(P+\h)^5(P+2\h)^5...(P+d\h)^5}  \]
is a formal vector-function with coefficients in this algebra. 

(b) Components of the vector-function form a fundamental solution  
to the linear differential equation (\ref{ODE}). }

\medskip

{\em Proof.} (a) Compute the coefficients of the series $L$ explicitly 
by the Duistermaat -- Heckman formula, 
\[ \int_{LX_d} \Phi (p,\h ) =\frac{1}{2\pi i} \oint \frac{\Phi (p,\h ) dp}
{p^5(p-\h)^5...(p-d\h)^5}, \] 
and change the order of summation. (b) Substitute $I$ into (\ref{ODE}). 
$\square $    

\medskip

{\bf The non-linear sigma-model.} The definition of the numbers $n_d$ accepted
in \cite{Gi1} was given by M. Kontsevich \cite{Kn} on the basis of intersection
theory in moduli spaces of stable maps to $X$. Let $X_{n,d}$ denote the moduli
orbifold of degree $d$ genus $0$ stable maps with $n$ marked points. For 
any such map $f:(\S,\e_1,...,\e_n)\to X$, the section $Q$ of the
degree $5$ line bundle over $X$ defines an element in the $5d+1$-dimensional
space $H^0(\S,f^*\calo(5))$ and hence --- a section of the orbibundle $V_d$
over $X_{n,d}$ formed by these spaces. According to Kontsevich's definition
the Euler class of this bundle capped with the fundamental class of the 
orbifold is taken on the role of the virtual fundamental class in the
corresponding space of stable maps to the quintic. Taking in account the
multiple cover formula \cite{AM} proved in \cite{M}, one gives the following
recursive definition of
the virtual numbers of degree $d$ spheres in quintic threefolds
\[ \sum_{m|d} \frac{n_{d/m}}{m^3} :=\int_{X_{0,d}} Euler (V_d) .\]
The first step in our proof of the conjecture extracting these numbers from 
the differential equation (\ref{ODE}) consists in mimicking Theorem $A$ 
within the framework of moduli spaces. 

Introduce the {\em graph space} $GX_{n,d}$ as the moduli orbifold of stable
maps to $GX=\CC P^1\times X$ of degree $d$ in projection to $X$ and degree
$1$ in projection to $\CC P^1$. Let $GV_d$ be the orbibundle with the
fiber $H^0(\S, f^*\ \calo_{\CC P^1}\otimes \calo_X(5))$ of dimension $5d+1$.
The space and the bundle inherit the $S^1$-action from the M\"obius
transformation group on $\CC P^1$. 

There is a natural equivariant birational isomorphism 
(see Main Lemma in \cite{Gi1})
\be \mu: GX_{0,d}\to LX_d  \label{mu} \end{equation}
defined as follows.

A stable map $f: \S \to \CC P^1\times X$ of bi-degree $(1,d)$ is represented
by the graph of a map $f_0: \CC P^1\to X$ of some degree $d_0\leq d$ and 
several ``vertical'' curves $f_i:(\S_i,\e_i) \to \{\z_i\}\times X$ of 
bi-degrees $(0, d_i)$, $i=1,...,r$, attached to the graph at the points
$f_i(\e_i)$. The number $r$ of such vertical curves
can vary from $0$ to $d$, and their total degree $d_1+...+d_r$ equals $d-d_0$.

The graph component $f_0$ is described by $5$ mutually prime binary forms
of degree $d_0$ uniquely up to a common constant factor  (we record them by 
$5$ polynomials $(x'_1(\z),...,x'_5(\z)$). Multiply these polynomials by the
same binary form of degree $d-d_0$ with the roots at $\z_1,...,\z_r $ of
multiplicity $d_1,...,d_r$ respectively (it is encoded by the polynomial
$x(\z)=(\z-\z_1)^{d_1}...(\z-\z_r)^{d_r}$ with the obvious convention about
the roots at $\z=\infty$). Then the polynomials $(x_1(\z),...,x_5(\z)):=
(x(\z)x'_1(\z),...,x(\z)x'_5(\z))$ represent the image $\mu (f)$ in the
space $LX_d$.

Obviously, the map $\mu $ is equivariant with respect to the natural actions of
the automorphism group $PSL_2(\CC )\times PSL_5(\CC )$ of $\CC P^1\times X$ on
$GX_{0,d}$ and $LX_d=Proj( H^0(\CC P^1, \calo_{\CC P^1} (d))\otimes 
H^0(X, \calo_X (1))^*)$ and is therefore independent on the choice
of the coordinate systems $(x_1:...:x_5)$ and $\z$ on $X$ and $\CC P^1$. 
In order to check that $\mu $ is regular let us assume for the moment that
the map $f$ is transverse to the five coordinate hyperplanes $\CC P^1\times 
x_i^{-1}(0)$, to $\infty \times X$, and that the intersections of $f(\S)$ with
the coordinate hyperplanes are away from $\infty \times X$. 
Then these intersections are all simple (even if $f$ is a multiple cover, they
are simple distinct non-singular points on $\S$),
their projections to $\CC P^1$ are away from $\infty$ and determine the roots 
of the polynomials
$x_1(\z),...,x_5(\z)$. The projective coordinates of the intersection point
$f(\S)\cap (\infty \times X)$ determine the top coefficients of the polynomials
$x_1(\z),...,x_5(\z)$ uniquely up to a non-zero scalar factor 
(we can fix it by requiring that the top coefficient of $x_1(\z)$ equals $1$).
Now, consider the point $[f]$ represented by the map $f$ in the moduli 
orbifold $GX_{0,d}$ and a local non-singular chart near this point. Due to the 
simplicity of the intersection points of $f(\S)$ with the coordinate 
hyperplanes and with $\infty \times X$ the roots and top coefficients of the
polynomials $x_i(\z)$ are regular functions of the map in a sufficiently small
neighborhood of $[f]$ (if this is not obvious yet, consider the universal
stable map $\ev: GX_{1,d}\to \CC P^1\times X$ and describe the intersection
points as local sections of the forgetting map $\ft: GX_{1,d}\to GX_{0,d}$). 
Thus the map $\mu $ is regular at $[f]$ since top coefficients and 
roots uniquely determine the polynomials. 

Choosing the coordinate systems on $\CC P^1$ and $X$ in general position 
to a given $f$ we conclude that $\mu $ is regular everywhere.

\medskip

The map $\mu $ allows one to compare the equivariant Euler classes of 
$LV_d$ and $GV_d$. Let $-p$ denote the equivariant 
Chern class of the Hopf bundle over $LX_d$ pulled-back to $GX_{0,d}$. 
Introduce the formal series 
\[ G(q,z):=
\sum_{d=0}^{\infty} q^d \int_{GX_{0,d}} e^{pz} Euler (GV_d).  \]

\medskip

{\bf Theorem B} {\em
\[ G(q,z)=\lan J(qe^{z\h},\h^{-1}), J(q,-\h^{-1})\ran \ ,  \]
where the formal vector-function $J$ with values in the cohomology algebra
$\QQ [P]/(P^4)$ is defined by
\be J(q,\h^{-1}):=e^{(P\ln q)/\h}[ 1+\h^{-1}\sum_{d=1}^{\infty}
q^d \ev_*(\frac{Euler (V_d)}{\h-c}) ], \label{J} \end{equation}
$c$ denotes the Chern class of the line orbibundle over $X_{1,d}$ formed
by cotangent lines to the curves at the marked point,
the map $\ev : X_{1,d}\to X$ is defined by evaluation of stable maps at the 
marked point, and the push forward is well-defined by 
\[ \lan \phi , \ev_* \Psi \ran = \int_{\CC P^4} 5P \phi \ev_*\Psi :=
\int_{X_{1,d}} (\ev^*\phi) \Psi  .\] }

\medskip

{\em Proof} (see \cite{Gi1}, Section $6$). It consists in application
of fixed point localization in equivariant cohomology. A fixed point of
the $S^1$-action on $GX_{0,d}$ is represented by a stable map $\S \to GX$
which consists of the graph of a constant map $\CC P^1\to X$ with two
``vertical '' curves of degrees $d-d'$ and $d'$ mapped to the slices $\z=0$ 
and $\z=\infty $ and attached to the graph. Thus the fixed point components 
are given by the diagonal constraint at the marked points in the Cartesian
product $X^{(0)}_{1,d-d'}\times X^{(\infty)}_{1,d'}$, $d'=1,...,d-1$, (and
are isomorphic to $X^{(0)}_{1,d}$ and $X^{(\infty)}_{1,d}$ for $d'=0,d$).
The normal bundle to the fixed point component has the equivariant Euler
class $(-\h)(-\h-c^{(0)})(\h )(\h-c^{(\infty)})$ 
(unless $d'=0$ or $d$ in which case a half
of the product should be taken). The Euler classes occur in the denominator
of the localization formula. The bundle $GV_d$ restricted to the fixed point
set splits into $V'_{d-d'}\oplus V'_{d'}\oplus \ev^*\calo (5)$ where $'$
symbolizes that the subspace in $H^0(\S, f^*\calo (5))$ consisting of
sections vanishing at the marked point is taken. Due to the
multiplicative property of Euler classes, $Euler (GV_d)$ factors 
correspondingly. Since the map $\mu $ sends the fixed point component in $GX_d$
to the corresponding fixed point component in $LX_d$, the class $p$ localizes
to $\ev^*(P)+d'\h$. It remains only to rearrange the summation over $d$ and 
$d'$ as the double sum over $d'$ and $d''=d-d'$. $\square $ 

\medskip

{\bf The divisor equation.} It is useful to figure out the place of the
series $J$ in the axiomatic structure of Gromov -- Witten 
theory (see for instance Dubrovin's book \cite{Db}). 
Genus $0$ Gromov -- Witten invariants of a K\"ahler manifold $Y$ define on 
$H^*(Y)$ the structure of a Frobenius manifold. It basically consists of
the quantum cup-product in each tangent space and a family of flat connections
$\nabla_{\h}$ on the tangent bundle defined by the quantum multiplication 
operators and
depending on a parameter which we denote $\h ^{-1}$. Structural constants of
the quantum cup-product are third partial derivatives $F_{\a\b\c}$ 
of the {\em potential}
\[ F(T)=\sum_{n,d} \frac{q^d}{n!}(T,...,T)_{n,d}\ ,  \]
where $T$ is a general cohomology class of $Y$,  
\[ (T,...,T)_{n,d}:=\int_{[Y_{n,d}]^{virt}} \ev_1^*(T)...\ev_n^*(T) ,\]
and the ``foreign'' formal parameter $q$ is introduced in order to provide
convergence of the $d$-sum in $l$-adic topology at least. 
For instance, the numbers $n_d$ for quintic threefolds are encoded by
the {\em Yukawa coupling}
\[ F_{ttt}=\sum_{n,d}\frac{q^d}{n!} (P,P,P,Pt,...,Pt)_{n+3,d} \]
which is obtained by differentiation of the potential restricted to
the $2$-nd cohomology space of the quintic and, as we will see soon, 
coincides with the formal Fourier series
\[ K(qe^t):=\sum_{d=0}^{\infty}  (P,P,P)_{3,d} q^de^{dt}=
5+\sum_{d=1}^{\infty} \frac{n_d d^3 q^d e^{dt}}{1-q^de^{dt}} .\]

The vector fields on $H^*(Y)$ flat with respect to the connections 
$\nabla_{\h}$ also represent some Gromov -- Witten invariants. 
Namely, the following matrix is the
fundamental solution to the system of linear differential equations 
$\nabla_{\h}S=0$:
\be  S_{\psi \phi}(T,\h^{-1}):=\lan \psi,\phi\ran + 
\sum_{n,d} \frac{q^d}{n!}(\psi,T,...,T,\frac{\phi}{\h-c})_{n+2,d} \ ,
\label{S} \end{equation}
where $\phi,\psi \in H^*(Y)$, and $c$ is the Chern class of the
universal cotangent line at the last marked point (as specified by the position
of this class in the correlator). 

Of course, all the correlators for quintic threefolds $Y$ are defined by means
of intersection theory in $X_{n,d}$ with the Euler classes of the bundles 
$V_d$ pulled-back from $X_{0,d}$ by forgetting maps.   

\medskip

{\bf Theorem C.} {\em
\be \lan J(qe^t,\h^{-1}), e^{-(P\ln q)/\h} \phi \ran =\lan 1, \phi \ran +
\sum_{n,d} \frac{q^d}{n!}(1, Pt,...,Pt, \frac{\phi }{\h -c})_{n+2,d} 
\label{desc} \end{equation}
} 

\medskip

{\em Proof} (\cite{Gi1}, Section $6$). It consists in application of the
{\em string equation} followed by an iterative application
of the {\em divisor equation}. The string and divisor equations read
respectively:
\be (1,T,...,T,\frac{\phi}{\h-c})_{n+2,d}=
\h^{-1}(T,...,T,\frac{\phi}{\h-c})_{n+1,d}, \label{string} \end{equation}
\be (P,T,...,T,\frac{\phi}{\h-c})_{n+2,d}=
d\ (T,...,T,\frac{\phi}{\h-c})_{n+1,d}
+ \h^{-1}(T,...,T,\frac{P\phi}{\h-c})_{n+1,d}, \label{divisor} \end{equation}
where one should require $n \geq 2$ in the case if $d=0$.
The equations allow to push the classes $\ev_i^*(1)$ and $\ev_i^*(P)$ of degree
$\leq 2$ from $X_{2+n,d}$ to $X_{1,d}$ along the forgetting maps
$X_{n+2,d}\to X_{n+1,d}\to ... \to X_{1,d}$. This transforms the sum of 
$d$-terms in the series (\ref{desc}) to 
$ \h^{-1} (\ \phi \exp(t(P+d\h)/\h)/(\h-c)\ )_{1,d} $
for $d>0$ (and to $\lan \exp(tP/\h),\phi \ran$ for $d=0$ since $X_{3,0}=X$).
$\square $

\medskip

A similar application of the divisor equation to the Yukawa coupling 
shows that $F_{ttt}(tP)=K(qe^t)$. On the other hand,
applying the string equation to the elements 
$S_{1,P}(tP,\h^{-1})$ of the matrix (\ref{S}) we extract the Yukawa coupling
as the coefficient at $\h^{-2}$ in the second derivative
in $t$. The theorem therefore implies that
\be K(e^t) = \lan \h^2 \frac{d^2}{dt^2} J (e^t,\h^{-1}), P\ran . 
\label{Yukawa} \end{equation}
Thus in order to prove the quintic formula it suffices to identify $J$ and
$I$ up to suitable mirror transformations.

\medskip

We complete this subsection by including a proof of the string and divisor
equations. The argument seems to be standard in Deligne -- Mumford theory.

\medskip

The maps $\ft : X_{k+1,d}\to X_k$ are defined \cite{Kn, BM} by forgetting 
the first marked point $\e_0$ of a stable map $f:(\S, \e_0,...,\e_k)\to X$ and
producing a new stable map $\tilde{f}: (\tilde{\S}, \e_1,...,\e_k)\to X$ by
contracting those components of $\S $ which have become unstable. 
The fiber of the forgetting map is canonically identified with the quotient
of the curve $(\tilde{\S},\e_1,...,\e_k)$ by the finite automorphism group
of the map $\tilde{f}$. In particular the map $\ft $ has $k$ canonical 
sections $\e_i$ 
defined by the marked points in the fibers and together with the
evaluation map $\ev_0 :X_{k+1,d}\to X$ can be considered as the universal
stable map with $k$ marked points.

One derives the equations (\ref{string}, \ref{divisor})
 by comparing the Chern class $c$ 
of the universal cotangent line bundle over $X_{1+k,d}$ with the pull-back
$\tilde{c}$ of the corresponding class from $X_{k,d}$. The cotangent lines
$T^*_{\e_k}{\S}$ and $T^*_{\e_k}\tilde{\S}$ are canonically identified 
{\em unless} $\e_k$ is situated on the same irreducible component $\S_0$ 
of $\S $ as $\e_0$, the map $f$ restricted to $\S_0$ is constant,
 $\S_0$ carries no other marked points and contains only one
singular point of the curve $\S$. Since $\S_0\simeq \CC P^1$ carries $3$
special points of different nature, the cotangent line $T^*_{\e_k}\S$
is canonically trivialized in this case. In terms of the universal cotangent
line bundles $l_k$ and $\tilde{l}_k=\ft^*l_k$ over $X_{k,d}$ with Chern 
classes $c$ and 
$\tilde{c}$ this means that $ \tilde{l}_k^*\otimes l_k$ has a canonical 
section non-vanishing
outside the divisor $\e_k(X_{k,d})$ defined by the universal marked point
$\e_k: X_{k,d}\to X_{k+1,d}$, and that $l_k$ restricted to this divisor is 
trivial. Since the restriction of $\tilde{l}_k$ to the divisor coincides
with the conormal bundle to this divisor by the very definition of the
universal cotangent line $T^*_{\e_k}\tilde{\S}$ on $X_{k,d}$, one concludes
that the section has the $1$-st order zero along the divisor. Thus
$\d=c-\tilde{c}$ is Poincare-dual in $X_{k+1,d}$ to the hypersurface 
$\e_k(X_{k,d})$, and $c \d  =0$.

We have
\[ (1,T,...,T,\tilde{c}^m\phi )_{k+1,d}=0,\  
(P,T,...,T,\tilde{c}^m\phi)_{k+1,d}=d(T,...,T,c^m\phi )_{k,d} \]
since the integral of the classes $1$ and $P$ over degree $d$ curves
equal $0$ and $d$ respectively. Finally, replacing $\tilde{c}$ by 
$c=\tilde{c}+\d$ yields the extra-terms
$(T,...,T,c^{m-1}\phi)_{k,d}$ and $(T,...,T,c^{m-1}P\phi)_{k,d}$.
This implies the string and divisor equations since $1/(\h-c)$ is the
eigenfunction of the operation $c^m\mapsto c^{m-1}$ with the eigenvalue
$1/\h$. $\ \ \square $.

\medskip

{\bf Torus action.} A link between $I$ and $J$ can be established 
via a recursion relation satisfied by their equivariant perturbations
$I^{eq}$ and $J^{eq}$. 

Consider the standard action of the $4$-dimensional torus $G$ on $X=\CC P^4$.
The equivariant cohomology algebra $H^*_G(X)$ is generated over the
coefficient ring $\QQ [\l]=H^*(BG)$ by the 
equivariant Chern class $-P$ of the Hopf line bundle and satisfies the
relation
\[ (P-\l_1)(P-\l_2)(P-\l_3)(P-\l_4)(P-\l_5)=0 \]
(where we assume that $\l_1+...+\l_5=0$).
Evaluation of an equivariant cohomology class on the invariant fundamental 
class can be computed via fixed point localization:
\[ \int_{[X]} \phi (P,\l)=\sum_{\a=1}^5 \frac{\phi (\l_{\a},\l)}{e_{\a}}=
\frac{1}{2\pi i}\oint \frac{\phi (P,\l )dP}{(P-\l_1)...(P-\l_5)} \]
where $e_{\a}=\Pi_{\b\neq \a}(\l_{\a}-\l_{\b})$ is the equivariant Euler
class of the tangent space to $X$ at the fixed point $r_{\a}$ where $P$ 
restricts to $\l_{\a}$. The equivariant Chern class of the anti-canonical 
line bundle $\calo_X(5)$ equals
$5P$, and we will denote $\lan \phi , \psi \ran $ the equivariant intersection 
pairing on $H^*_G(X)$ with this Chern class:
\[ \lan \phi , \psi \ran := \int_{[X]} 5P \phi \psi .\]

The torus $G$ acts on the moduli orbifolds $X_{n,d}$ which allows one to 
introduce the equivariant Gromov -- Witten correlators (such as 
$(T,...,T)_{n,d}$ and $(\psi, T,..., T, \phi c^k)_{n+2,d}$
where $c$ and $\phi, \psi, T\in H^*_G(X)$ are equivariant cohomology classes. 
The correlators take values in $\QQ [\l ]$ and turn into corresponding
non-equivariant Gromov -- Witten invariants when specialized to $\l=0$.

Similarly, one can use $S^1\times G$-equivariant intersection theory in
the spaces $LX_d$ and $GX_{0,d}$ and carry over Theorems $A(a),B,C$ to the 
equivariant setting. The proofs are identical to those given above, but
the equivariant counterpart of the hypergeometric
series $I$ is defined now by the series
\be I^{eq}:=e^{(P\ln q)/\h}\sum_{d=0}^{\infty} q^d 
\frac{\Pi_{m=1}^{5d} (5P+m\h)}
{\Pi_{m=1}^d \Pi_{\a=1}^5 (P-\l_{\a}+m\h)}, \label{I} \end{equation}
with coefficients in the {\em equivariant} cohomology algebra of $X$  
generated by $P$ (and represents a fundamental solution to some $5$-th order
ODE, which is however irrelevant for our goal in this paper).

In order to describe the recursion relation satisfied by the hypergeometric
series $I^{eq}$ let
us strip off the factor $\exp (P\ln q)/\h $, denote the remaining series by 
$Z^{(hg)}$ 
and denote by $Z^{(hg)}_{\a}$ its fixed point localizations:
\[ Z^{(hg)}_{\a} (q,\h^{-1},\l )= \sum_{d=0}^{\infty}
\frac{q^d}{d!\h^d}\frac{\Pi_{m=1}^{5d} (5\l_{\a}+m\h)}
{\Pi_{m=1}^d \Pi_{\b\neq\a} (\l_{\a}-\l_{\b}+m\h)} \ .\]
Coefficients of the formal $q$-series $Z^{(hg)}_{\a}$ are degree $0$ 
rational functions in $\h $ with the first order pole at 
$\h=(\l_{\b}-\l_{\a})/m$, $\b\neq \a$,
$m=1,...,d$, and a high order pole at $\h=0$. Rewriting the rational functions
as sums of elementary fractions we arrive at the following recursion relation:
\be
Z_{\a}(q,\h^{-1},\l)= \label{rec} \end{equation}
\[ 1+\sum_{d=1}^{\infty} q^d\frac{R_{\a,d}(\h, \l )}{\h^d}
+  \sum_{\b\neq\a}\sum_{m=1}^{\infty} C_{\a}^{\b}(m) 
\frac{q^m}{\l_{\a}-\l_{\b}+m\h} Z_{\b}(q, \frac{m}{\l_{\b}-\l_{\a}},\l), \]
where $R_{\a,d}$ are some degree $\leq d$ polynomials in $\h$ with coefficients
rational in $\l$, and $C_{\a}^{\b}(m)$ are rational functions in $\l$ which
we will call the {\em recursion coefficients}. Since 
$Z^{(hg)}_{\b}\equiv 1 (\mod q)$, 
the coefficient $C_{\a}^{\b}(d)$ can be read off
the $q^d$-term of the series $Z^{(hg)}_{\a}$ as the residue at the pole 
$ \h = (\l_{\b}-\l_{\a})/d$.

We will refer to the sequence of polynomials $R_{\a,d}$ as the {\em initial
condition}: given such an initial condition, the recursion relation allows
to recover the solution $\{ Z_{\a}, \a=1,...,5 \}$ unambiguously.

\medskip

Now, starting with the Gromov -- Witten invariant $J^{eq}$, 
introduce the vector $q$-series 
$Z^{(GW)} (q, \h^{-1},\l)$ with coefficients in the algebra $H^*_G(X,\QQ (\h))$
by stripping off the factor $\exp (P\ln q)/\h $,  and denote 
by $Z^{(GW)}_{\a}$ the localization of $Z^{(GW)}$ at the fixed $r_{\a}$.

\medskip

{\bf Theorem D} (\cite{Gi1}, Section $11$). {\em
The series $\{ Z^{(GW)}_{\a}, \a=1,...,5 \}$ satisfy the recursion relation
(\ref{rec}) with the same recursion coefficients $C_{\a}^{\b}(d)$ (and with
another initial condition).}  

\medskip

{\em Proof} ( \cite{Gi1}, Sections $9,11$). It is based on localization
to fixed points of the torus $G$ action on the moduli orbifolds $X_{1,d}$.
Consider a stable map $f: (\S, \e) \to X$ representing 
such a fixed point. The combinatorial structure of the curve $\S $ is
described by a tree of irreducible components (isomorphic to $\CC P^1$ each).
Some components are mapped onto the straight lines in $X=\CC P^4$ connecting
the fixed points $r_{\a}, \a=1,...,5 $, of the torus action, 
and the map is a multiple cover $\z \mapsto \z^m$ in suitable 
affine coordinates on the source and
target $\CC P^1$, so that $\z=0, \infty$ are mapped to the fixed points.
The remaining irreducible components of $\S $ 
are mapped to the fixed points in $X$. The marked point $\e $ must be
mapped to one of the fixed points $r_{\a}$.

The fixed point in $X_{1,d}$ represented by $f$ does not contribute to
$Z^{(GW)}_{\a}$ via localization formulas {\em unless} 
$f(\e)=r_{\a}$.    

Suppose that $\e$ is situated in an irreducible component of $\S$
mapped to $r_{\a}$. Consider the whole connected component of the fixed point
set $X_{1,d}^G$ in $X_{1,d}$ which contains the equivalence class $[f]$. 
We will show
that this connected component contributes to $Z^{(GW)}_{\a}$ by a polynomial in
$\h^{-1}$. Indeed, the component can be described as the (quotient by
a finite group of the) product of some Deligne -- Mumford spaces 
$\M_{0, k}$ (why? --- see \cite{Kn} where the fixed point set is described). 
The universal cotangent line orbibundle over $X_{1,d}$
restricted to the connected component of $X_{1,d}^G$
coincides with the universal cotangent line at one of the marked points over
one of the factors $\M_{0,k}$. Thus its Chern class $c$ is nilpotent on
this component, and the geometrical series $(\h-c)^{-1}$ reduces to a finite
sum of terms $c^l/\h^{l+1}$. Notice that $l\leq \dim \M_{0,k}=k-3$ is 
bounded by the total degree of the map $f$.  

The fixed point localization terms just discussed form the initial condition
in (\ref{rec}). We will show that contributions of all other fixed points can
be arranged as the recursive part of (\ref{rec}). The idea is to cut off the
component of the curve $(\S ,\e) \to X$ carrying the marked point and to 
observe that the rest of the curve represents a torus-invariant curve of 
smaller degree.    

In greater detail, suppose that $\e $ is situated at $\z=0$ on a multiple 
cover $\z \mapsto \z^m$
of the line connecting $r_{\a}$ with $r_{\b}$. Then the universal cotangent
line orbibundle restricted to the connected component of $X_{1,d}^G$ is 
topologically trivial (since $(T^*_{\e}\S)^{\otimes m}$ coincides in this case
with the cotangent space $T^*_{r_{\a}}\CC P^1$ to the line joining $r_{\a}$
and $r_{\b}$), but it carries a nontrivial infinitesimal action of $G$ given
by the character $(\l_{\a}-\l_{\b})/m$. Thus the localization of $(\h -c)^{-1}$
at this fixed point component yields the simple fraction 
$m q^m (m\h + \l_{\b}-\l_{\a})^{-1}$. The factor $m$ is eventually 
compensated by
the order of the automorphism group of the map $\z\to \z^m$ which occurs
in the denominator of localization formulas on orbifolds. The weight $q^m$
counts the degree of this map as a curve in $X$.   

The whole contribution of the fixed point component to $Z^{(GW)}_{\a}$ via 
localization formulas includes two more factors. Each of them takes in account
the equivariant Euler classes of the orbibundle $V'_d$ 
and of the normal orbibundle
to $X_{1,d}^G$ which occurs in the denominator of localization formulas. 

The first factor corresponds to the irreducible component $C=\CC P^1$
of $\S $ carrying the marked point, and the second one --- corresponds to
the remaining part $\tilde{\S}$ of the curve $\S $. The map $f$ restricted to
$\tilde{\S}$ has degree $d-m$ and represents a point in the space
$X_{1,d-m}^G$. The fiber $H^0(\S , f^*V')$ of $V'_d$ contains the subspace
of sections vanishing on $C$, which coincides with the fiber of $V'_{d-m}$.
The normal spaces to $X^G_{1,d}$ split similarly into parts corresponding to
$C$ and $\tilde{\S}$.
The intersection point $\z=\infty $ of $C$ with $\tilde{\S}$ 
plays the role of the marked point $\tilde{\e}$ in $\tilde{\S}$. 
The deformation of $f$ corresponding to smoothening of the curve $\S $ at 
the double point $\tilde{\e}$ is represented in the tangent space to $X_{1,d}$
by the line $T_{\tilde{\e}}C\otimes T_{\tilde{\e}}\tilde{\S}$ and contributes
the factor $(\l_{\b}-\l_{\a}) m^{-1}-\tilde{c}$ to the denominator of the 
localization formula. Thus the contribution of $\tilde{\S}$ is correctly
accounted by the factor $Z^{(GW)}_{\b}(q,m/(\l_{\b}-\l_{\a}),\l)$ 
in the recursion relation (\ref{rec}).
    
The remaining factor in the localization formula is a rational function
of $\l $ and can be computed explicitly as the ratio of two Euler classes
in the case $m=d$ when $\tilde{\S}$ is a point (in the example of quintics 
we are studying, the factor has actually been computed in \cite{Kn}). 
It turns out to 
coincide with the recursion coefficient $C_{\a}^{\b}(m)$.  
We recommend the reader to carry out this computation (it amounts to 
analyzing the torus action on spaces of holomorphic sections of 
$\calo (5)$ and $T_X$ lifted to the multiple cover $\z \mapsto \z^m$
of the line joining $r_{\a}$ and $r_{\b}$) or at least to look at some
details of this computation in \cite{Kn}. $\square $

A plausible argument in \cite{Gi1, Gi3} intended to explain the ``miraculous''
coincidence of the recursion coefficients has been formalized in \cite{LLY}.

\medskip

{\bf Polynomiality.} The recursion relation (\ref{rec}) has much more
solutions with various initial conditions than the mirror transformations
can handle. However, according to Theorem A and Theorem B, the solutions
$Z^{(hg)}$ and $Z^{(GW)}$ have the following polynomiality property.

Let us call a solution $Z (q,\h^{-1},\l)$ to the recursion
relation (\ref{rec}) {\em polynomial} if the formal $(q,z)$-series
\be \lan  Z(qe^{\h z}, \h^{-1}, \l ), e^{Pz} \ Z(q, -\h^{-1},\l ) \ran 
\label{pol} \end{equation}
has coefficients polynomial in $\h $.

\medskip

The solution $Z^{(GW)}$, by the very definition 
(\ref{J}), satisfies also the asymptotical condition
\be Z(q,\h^{-1},\l ) = 1 + o (\h^{-1}). \label{as} \end{equation}    

\medskip

{\bf Theorem E} (\cite{Gi1}, Proposition $11.5$).  
{\em A polynomial solution to the 
recursion relation (\ref{rec}) satisfying the asymptotical condition 
(\ref{as}) is unique.}

\medskip

{\em Proof:} perturbation theory. Let $Z$ be a polynomial solution 
and let $\d R=R_d-R^{(GW)}_d$ denote
the discrepancy in the initial conditions for $Z$ and $Z^{(GW)}$ 
with minimal $d>0$. Then $Z$ and $Z^{(GW)}$ coincide modulo $q^d$ due to
the recursion relation. The polynomiality property for $Z$ and $Z^{(GW)}$
modulo $q^{d+1}$ translates into regularity at $\h=0$ of 
\[ \lan \d R (\h)\h^{-d}, e^{(P+d\h)z} \ran + \lan \d R(-\h)(-\h)^{-d},
e^{Pz}\ran .\]
Localizations of this intersection index to fixed points in $X$ are ---
for each power of $\h^{-1}$ --- finite sums of monomials 
$z^l\exp (\l_{\a} z),\  l=0,1,2,...,\ \a=1,...,5 $. Linear independence
of such monomials for generic $\l $ implies that the localizations
\[ \d R_{\a} (\h )\h^{-d} e^{\h z} + \d R_{\a}(-\h)(-h)^{-d}, \]
where $R_{\a}(\h)$ are some polynomials in $\h $ of degree $\leq d$,
must be regular at $\h =0$ on their own. 

Consider first $(A\h^{-2}+B\h^{-3})\exp (\h z )+(A\h^{-2}-B\h^{-3})$ modulo
$z^2$. The regularity condition at $\h =0$ implies $A=0$ and then $B=0$.
Applying this argument inductively we conclude that $\d R_{\a}(\h) \h^{-d}=
A_{\a}+B_{\a} \h^{-1}$ where $A_{\a}, B_{\a}$ do not depend on $\h$. 

Assuming now that $Z$ also satisfies the asymptotical condition (\ref{as}) we
find $\d R=0$.   $\square $

\medskip

{\bf Mirror transformations.} The hypergeometric series $I^{eq}$ 
has the asymptotical expansion
\[ I^{eq}=e^{(P\ln q)/\h }(f_0(q)+f_1(q)\frac{P}{\h }+o(\h^{-1})) ,\]
where the series 
\[ f_0=\sum_{d=0} q^d \frac{(5d)!}{(d!)^5} ,\]
\[ f_1=\sum_{d=1} q^d \frac{(5d)!}{(d!)^5}(\sum_{m=d+1}^{5d} \frac{5}{m}) ,\]
are found from (\ref{I}) (remember that $\l_1+...+\l_5=0$).

The {\em mirror transformations}, namely the division of $I^{eq}$ by $f_0$
followed by the change of variable $\ln q \mapsto \ln q \ +\ f_1(q)/f_0(q) $,
transform $I^{eq}$ to a new vector-function with the same asymptotical 
behavior as 
\[ J^{eq} = e^{(P\ln q)/\h} (1+o(\h^{-1})) .\]
Therefore 
the following theorem guarantees that the transformed series coincides
with $J^{eq}$. Passing to the non-equivariant limit $\l=0$ we conclude that 
the same mirror transformations take $I$ into $J$. 

Thus the Yukawa coupling 
(\ref{Yukawa}) is indeed extracted from the fundamental solution $I$ to the 
differential equation (\ref{ODE}) by the procedure conjectured in \cite{COGP}.
 
\medskip

{\bf Theorem F} (\cite{Gi1}, Propositions $11.3, 11.6$). {\em The mirror 
transformations take polynomial solutions of the recursion relation (\ref{rec})
to polynomial solutions of the same recursion relation.}  

\medskip

{\em Proof:} straightforward (see \cite{Gi1}). 
The division operation does not change the
form of the recursion relation and also preserves the polynomiality property
since the extra factor $f_0(q\exp (\h z)) f_0(q)$ does not produce negative
powers of $\h $ in the $(q,z)$-series (\ref{pol}). The change 
of the variables $\ln q \mapsto \ln q + g(q)$  transforms $z$ in this
series into $z+ [g(q\exp (\h z)) - g(q)]/\h $. At $\h=0$ the difference 
vanishes. It is therefore divisible by $\h$, and thus the change of variables
preserves the polynomiality property too. 

When applied to the
recursion relation (\ref{rec}) literally, the change of variables modifies 
it to a new recursion relation. In the new form the elementary
fraction 
\[ \frac{q^m}{m \h +\l_{\a}-\l_{\b}} = \frac{q^m}{\h (\l_{\a}-\l_{\b})}
\ \ \frac{1}{\h^{-1}+m/(\l_{\a}-\l_{\b})} \]
occurs with the extra factor
\[ \d :=
\exp [ m g(q)+\l_{\a} g(q)\h^{-1} - \l_{\b} g(q) \frac{m}{\l_{\b}-\l_{\a}}]\]
\[= \exp [ \l_{\a} g ( \h^{-1}+m/(\l_{a}-\l_{\b}) ) ] . \]
Thus $\d -1$ is divisible by $\h^{-1}+m/(\l_{\a}-\l_{\b})$. Since $g(0)=0$,
the result of this division is a $q$-series with coefficients polynomial in
$\h^{-1}$ at each power of $q$.    
Thus the transformation  affects only the initial condition and takes 
a solution of the recursion relation into another solution. $\square $

\bigskip

{\bf Updates.} 

\medskip

{\em Definitions.} The definition of virtual numbers $n_d$ in terms of 
Euler classes of the orbibundles $V_d$ over $X_{0,d}$ should be considered
as tentative and has been replaced by a more universal construction, due to
J. Li \& G. Tian \cite{LT}, 
of virtual fundamental classes $[Y_{0,d}]$ defined
in intrinsic terms of the quintic $3$-fold $Y$ rather than in terms of the
embedding $i: Y\subset X$. 
Thus in order to place the above proof of the
quintic formula into the framework of contemporary definitions one needs to 
check that $i_*[Y_{0,d}]$ in $H_*(X_{0,d})$ equals the cap-product of the 
virtual fundamental class $[X_{0,d}]$ with $Euler(V_d)$. This is easy and can
be done as follows.

Consider first a model problem studied in \cite{F}: given a holomorphic
section $s: B\to E$ of a vector bundle $E\to B$ over a compact complex 
manifold, construct a cycle
in the zero locus $Z=s^{-1}(0)$ 
of the section which is Poincare-dual to the Euler class of
the bundle. The model problem is solved by the normal cone construction: the 
normal cone $C \subset E |_Z$ to the zero locus has pure dimension 
$\dim B$, and the homological intersection in $E |_Z$ 
of the fundamental cycle $[C]$ with that of the zero section represents in
$H_*(Z)$ the required Euler class. The construction can be adjusted to
the orbifold/orbibundle situation.

The virtual fundamental cycle construction in \cite{LT} is based on the 
observation that the normal cone $C \subset E |_Z$ is intrinsic with respect
to the scheme structure of $Z$ and the {\em tangent-obstruction} complex 
$ds: T_B|_Z \to E|_Z$ of vector 
bundles over $Z$ defined by the differential of the section. The kernel of
$ds|_{z\in Z}$ is the algebraic tangent space $T_zZ$. In the case when $Z$
is a moduli space of stable maps to $Y$, J. Li \& G. Tian exhibit a 
tangent-obstruction complex
of orbibundles $T\to E$ with this property and by this define the intrinsic 
normal orbi-cone $C\subset E$ and the virtual fundamental class 
$[Z]:= [C]\cap \ [zero\ section]$ in $H_*(Z,\QQ)$. 

In our situation $Z=Y_{0,d}$ is given in the orbifold $X_{0,d}$ by a section 
$s$ of $V_d: E_{0,d}\to X_{0,d}$. Using the exact sequence 
$0\to T_Y \to T_X|_Y \to N_Y\to 0$ where $N_Y$ is the normal bundle to $Y$ in 
$X$ we obtain the exact sequence 
\[ 0\to H^0(\S, f^*T_Y)\to H^0(\S, f^*T_X)\to
H^0(\S, f^*N_Y)\to H^1(\S, f^*T_Y)\to 0 \]
for each stable map $f:\S \to Y$. Since the fibers of $V_d|_Z$ coincide with
$H^0(\S, f^*N_Y)$, this implies (via the description \cite{Kn, LT} of the 
algebraic tangent space $T_{[f]} Y_{0,d}$ in cohomological terms of 
deformation theory) that the complex 
$ds: T_{X_{0,d}}|_Z \to E_{0,d}|_Z$ can be taken on the role of the 
tangent-obstruction complex in the definition of the virtual fundamental 
cycle $[Y_{0,d}]$. Thus
this cycle represents in $H_*(X_{0,d})$ the Euler class of $V_d$.

While this obvious argument shows that the GW-invariant of $Y$ in question
can be computed in terms of the GW-theory for the convex bundle $V: E\to X$,
some other GW-invariants of $Y$ can not be interpreted in terms of the bundle.
In order to distinguish the GW-theory of the bundle from the GW-theory of $Y$
we, following A. Schwarz, refer in \cite{Gi3} to the first one as the
GW-theory of the {\em supermanifold} $\Pi E$. 

\medskip

{\em The map $\mu $.} Applications of our approach to complete 
intersections in toric manifolds more general than projective spaces showed
that some steps in the above proof are redundant. 
In the remaining part of the text we discuss several such steps which can
be simplified or avoided. 
The first of them is our use of the map $\mu : GX_{0,d}\to LX_d$.

The map was used in order to define the equivariant class $p$ on $GX_{0,d}$ 
as a pill-back of the corresponding class on $LX_d$ and thus assure the
polynomiality property of the GW-invariant $J^{eq}$. Consider instead the
following $S^1\times G$-equivariant GW-invariant of $\CC P^1\times X$:
\be \sum_{d=0}^{\infty} \sum_{n=0}^{\infty} \frac{q^d}{n!} 
[ T,...,T ]_{n,d}, \label{n/d} \end{equation}
where $T=z (p\otimes P)$, $p$ is the generator of the $S^1$-equivariant 
cohomology of $\CC P^1$ satisfying $p(p-\h)=0$, $P$ is the generator of
the $G$-equivariant algebra of $X$, and $[ ... ]_{n,d}$ is the equivariant
GW-invariant defined by integration over $GX_{n,d}$ against the Euler class
of $GV_{n,d}$.
The series (\ref{n/d}) is defined without fixed point localization and
thus is a $(q,z)$-series with coefficients polynomial in $\h$ and $\l $.
On the other hand, applying localization to fixed points of $S^1$-action on
$GX_{n,d}$ as in the proof of Theorem B (notice that $p$ localizes to
$0$ at $\z=0$ and to $\h $ at $\z=\infty $) and then using the divisor 
equation for $zP$ as in the proof of Theorem $C$ we will find that the series
(\ref{n/d}) coincides with 
$\lan J^{eq}(q\exp \h z, \h^{-1}), J^{eq}(q,-\h^{-1}) \ran $.
This argument was mentioned in \cite{Gi3} and was used in \cite{Gi4}.

\medskip

{\em Theorem F.} The invariance of the recursion
relation  (\ref{rec}) under mirror transformations can be deduced from the
string and divisor equations. Namely, consider the $G$-equivariant 
GW-invariant $\calJ^{eq}$ defined by integration over $X_{n+2,d}$ against the
Euler classes of $V_{n+2,d}$:
\be \lan \calJ , \phi \ran = \lan 1, \phi \ran +\sum_{(n,d)\neq (0,0)}
\frac{q^d}{n!} (1, T, ..., T, \frac{\phi e^{P\ln q/\h}}{\h - c})_{n+2,d}\  
\label{calJ} \end{equation}
with  $T=a(q)+b(q)P$ where $a,b$ are power $q$-series 
vanishing at $q=0$. One can derive a recursion relation for $\calJ^{eq}$ in 
exactly the same way as we derived the recursion relation for $J^{eq}$ in 
Theorem D. The recursion coefficients will be the same as in Theorem D,
but the initial condition will depend now on $a$ and $b$. On the other hand 
the string and divisor equations show that 
\[ \calJ^{eq}=e^{a(q)} J^{eq} (qe^{b(q)},\h^{-1}) \]
and is therefore a result of a mirror transformation applied to $J^{eq}$.
Also, using the
argument from the previous subsection with suitable function of $a$ and $b$ 
on the role of $T$, we can make sure that $\calJ^{eq}$ must {\em a priori} 
satisfy the polynomiality condition as well.  

\medskip

{\em Mirror transformations.} As it was shown in \cite{Gi1} 
(Section $12$) and \cite{Gi4} (Section 5), there is a ``non-linear Serre
duality'' equivalence between genus $0$ equivariant GW-theory for a convex 
supermanifold $\Pi E$ and such a theory for the non-compact total space 
$E^*$ of the concave dual bundle  $V^*: E^*\to X$. Namely, their genus $0$ 
GW-invariants differ by a change of variables which can be explicitly described
in terms of GW-invariants of either of them (see \cite{Gi4}, Section 5).

On the other hand, the mirror formulas can be generalized, as it was shown in
\cite{LLY} and \cite{Gi4} (Section $4$), to include genus $0$ equivariant 
GW-invariants of concave bundles $E^*$. The proof is completely parallel to
the one given above. However, in the case if the bundle $V^*$ is the direct 
sum of at least two line bundles, the GW-invariant $J^{eq}_{E^*}$ for $E^*$ is
{\em equal} to the corresponding hypergeometric series $I^{eq}_{E^*}$ which in
this case happens to satisfy the asymptotical condition of the uniqueness 
Theorem $E$. With this observation, a proof of the quintic formula looks
as follows (see Section $5$ in \cite{Gi4}).

Describe quintic $3$-folds by two equations in $\CC P^5$ of degree $1$ and $5$.
For the concave bundle $\calo(-1)\oplus \calo(-5)$ over $\CC P^5$, prove the
equality $J^{eq}_{E^*}=I^{eq}_{E^*}$ following the steps 
Theorem B --- Theorem D --- Theorem E as explained above. Now the mirror
transformation between $J^{eq}_{\Pi E}$ and $I^{eq}_{\Pi E}$ emerge from
the general formulas of ``non-linear Serre duality''.

\medskip

{\em Theorem D.} The proof of mirror formulas given in \cite{LLY}
is based on some recursive property of the $S^1\times G$-equivariant Euler 
classes of the bundles $GV_d$ over $GX_{0,d}$ named {\em eulerity}. In fact 
the property can be easily deduced from our recursion relation (\ref{rec})
as it is done in \cite{Gi1} (Proposition $11.4$ (b)). The inverse implication
is also immediate. However, the proof of
the eulerity property given in \cite{LLY} is based on an argument similar to
our proof of Theorem $B$ (localization for the $S^1$-action) and thus 
completely eliminates the role of localization formulas for the $G$-action 
as a computational tool (and uses only the very fact of their existence).

With this observation, the reduction of the quintic formula to ``non-linear 
Serre duality'' theorem looks particularly short: prove eulerity property
in the case of the concave bundle $\calo(-1)\oplus \calo(-5)$ over $\CC P^5$ 
following the argument in \cite{LLY} and use a uniqueness theorem (parallel to
Theorem $E$ above) in order to identify $J^{eq}_{E^*}$ with $I^{eq}_{E^*}$.

\bigskip

We see that the only
new, after M. Kontsevich's paper \cite{Kn}, geometrical construction which has
survived so far through all variants of the proof of the quintic formula is
the {\em $S^1$-equivariant} theory on the graph spaces $GX_{0,d}$ which 
originates from the loop space interpretation \cite{Gi2} (see also \cite{V})
of Gromov -- Witten invariants.

All other ingredients of the proof have somewhat combinatorial character and 
can be interchanged and simplified. 
This progress does not mean however that the 
mirror symmetry phenomenon has been adequately understood.

\enddocument